\documentclass[review]{elsarticle}

\usepackage{hyperref}
\usepackage{amssymb,amsmath,amsthm}
\usepackage{color}
\usepackage{graphicx}
\usepackage{subfigure}
\usepackage{epstopdf}
\usepackage{graphics,epsfig}

\newcommand{\argmin}{\arg\!\min}
\renewcommand{\qedsymbol}{$\blacksquare$}

\journal{}

\begin{document}

\begin{frontmatter}

\title{Distributed Real-Time Non-Linear Receding Horizon Control Methodology for Multi-Agent Consensus Problems}

\author[mymainaddress]{Fei Sun\corref{mycorrespondingauthor}}
\author[mymainaddress]{Kamran Turkoglu}

%
\cortext[mycorrespondingauthor]{Email: fei.sun@sjsu.edu}
%
\address[mymainaddress]{Aerospace Engineering, San Jos\'e State University, San Jose, CA 95192}

\begin{abstract}
This work investigates the consensus problem for multi-agent nonlinear systems through the distributed real-time nonlinear receding horizon control methodology. With this work, we develop a  scheme to reach the consensus for nonlinear multi agent systems under fixed directed/undirected graph(s) without the need of any linearization techniques. For this purpose, the problem of consensus is converted into an optimization problem and is directly solved by the backwards sweep Riccati method to generate the control protocol which results in a non-iterative algorithm. Stability analysis is conducted to provide convergence guarantees of proposed scheme. In addition, an extension to the leader-following consensus of nonlinear multi-agent systems is presented. Several examples are provided to validate and demonstrate the effectiveness of the presented scheme and the corresponding theoretical results.
\end{abstract}


\end{frontmatter}


\section{INTRODUCTION}

With their sophisticated structure, multi-agent related consensus problems have attracted significant interest in recent years. The complex nature and sophisticated framework of multi agent consensus problem serves as a fertile ground for the application of advanced control algorithms, and found basis in many areas including cooperative control, formation control, flocking and rendezvous \cite{saber2004,dunbar2006,dunbar2007,hui2008,keviczky2008,you2011,yu2011,karl2012,xie2014,WangChenMa2014,wang2015,zhaodw2015,zhao2015}. 

So far, the consensus methodologies have been widely explored for multi-agent \emph{linear} dynamical systems \cite{LiDuanChenHuang2010,ZhuYuan2014,MovricLewis2014}. However, most physical systems are nonlinear in their nature or even exhibit the phenomenon of chaos. The combination of consensus and nonlinear dynamics remains as a major challenge in literature, although some preliminary results have been presented in recent years. For instance, in the work of \cite{BaussoGiarrePesenti2006}, a new rule was introduced for nonlinear protocol design allowing consensus on a general set of values. In \cite{QuChunyuWang2007}, the investigation was concentrated on nonlinear cooperative control for consensus of nonlinear and heterogeneous systems. \cite{LiuXieRenWang2013} discussed the problem of consensus for multi-agent systems with inherent nonlinear dynamics under directed topologies, where a variable transformation method is used to convert the consensus problem to a partial stability problem. \cite{XuCaoYuLu2013} studied the leader-following consensus of nonlinear multi-agent systems without the assumption that the topology is connected or fixed.

On the other hand, Receding Horizon Control based methodologies gained significant momentum in the last two decades which are used to obtain an optimal feedback control law by minimizing the desired performance index for a given finite horizon. In that sense, receding horizon control is one of the powerful methodologies that is adapted to consensus problem of multi-agent dynamics. Based on this approach, there have been many results developed for consensus problems and its applications. The work of \cite{dunbar2005} presented a distributed receding horizon control law for dynamically coupled nonlinear systems based on its linearization representative. The robust distributed receding horizon control methods were studied in \cite{LiShi2013} for nonlinear systems with coupled constraints and communication delays. \cite{shi2014} proposed a robust distributed model predictive control methods for nonlinear systems subject to external disturbances. In a very recent study \cite{QiuDuan2014}, brain-storm type of optimization was combined with receding horizon control strategies for UAV formation flight dynamics.

Most of the previously conducted studies on nonlinear consensus problem either employ some sort of adaptation mechanisms (through adaptive control methodologies) or simplistic linearization methods. In this paper, we propose a scheme for solving nonlinear consensus problem by utilizing the real-time nonlinear receding horizon control (NRHC) methodology, which avoids such pitfall and utilizes complete architecture of the nonlinear system. For a given fixed directed/undirected network, the nonlinear consensus problem is converted to an optimization framework and the control protocol is designed locally for each agent by real-time nonlinear receding horizon control algorithm which is an important extension to the existing literature and presents the novelty of proposed NRHC procedure. Here, a non-iterative optimization algorithm is applied to avoid high computational complexity and large data storage. During this process, each agent only needs to obtain its neighbor's state \emph{once} via the given network at each time instant, which is more efficient than the other distributed strategies that involve multiple information exchanges and predicted trajectories of states. This presents the second major novelty proposed with this approach. With this formulation, we also provide asymptotic stability guarantees to achieve consensus within the given topology. The proposed scheme is then extended to the leader-following consensus problem of nonlinear multi-agent systems. The effective nature of the proposed methodology is demonstrated through several examples with different topologies. With this approach, we demonstrate the applicability of the distributed real-time nonlinear receding horizon control as a consensus routine on the nonlinear multi-agent systems.

The organization of the paper is as following: In Section-\ref{sec:prob_formulation}, some preliminaries and problem statement of the distributed real-time nonlinear receding horizon control is presented for nonlinear consensus dynamics. The framework and the detailed control design are investigated in Section-\ref{sec:dist_control}. Stability analysis is discussed in Section-\ref{sec:stab_analysis}. In Section-\ref{sec:simulation}, theoretical aspects of the proposed methodology are applied on Lorenz oscillators and L\"{u} oscillators examples. And with Section-\ref{sec:conclusions}, the paper is concluded.

\section{PRELIMINARIES AND PROBLEM STATEMENT}\label{sec:prob_formulation}

In this section, we introduce the notation used throughout this paper and we provide a broad definition of the main problem studied. Later on, we present algorithms for solving this problem.

In this work, $R$ denotes the real space. For a real matrix $A$, its transpose and inverse are denoted as $A^T$ and $A^{-1}$, respectively. The symbol $\otimes$ represents the Kronecker product. For matrices $X$ and $P$, the Euclidean norm of $X$ is denoted by $\|X\|$ and the $P$-weighted norm of $X$ is denoted by $\|X\|_P=\sqrt[P]{X^TPX}$. $I_n$ stands for the identity matrix of dimension $n$. Given a matrix $P$, $P>0$ ($P<0$) represents that the matrix is positive definite (or negative definite). Here, we define the column operation $\text{col}(x_1,x_2,\cdots,x_n)$ as $(x^T_1,x^T_2,\cdots,x^T_n)^T$ where $x_1,x_2,\cdots,x_n$ are column vectors.

Consider a multi-agent system of $M$ nonlinear agents. For each agent $i$, the dynamic system is given by:
\begin{equation}\label{model}
\dot{\mathbf{x}}_i=f(\mathbf{x}_i,\mathbf{x}_{-i},t)=\mathbf{F}(x_i)+\mathbf{u}_i,
\end{equation}
where $\mathbf{x}_{i}=(x_{i1},x_{i2},\cdots,x_{in})^{T}$ is the
state vector of the $i$th oscillator, $\mathbf{x}_{-i}$ are the collection of agent $i$'s neighbor's states, the function $\mathbf{F}(\cdot)$ is the corresponding nonlinear vector field, and $\mathbf{u}_i$ is the control input of agent $i$. Here, function $\mathbf{F}(\cdot)$ satisfies the global Lipschitz condition. Therefore there exists positive constant $\beta_i$ such that
\begin{equation*}
\|\mathbf{F}(\mathbf{x}_i)-\mathbf{F}(\mathbf{x}_j)\|\leq \beta_i\|\mathbf{x}_i-\mathbf{x}_j\|.
\end{equation*}
This condition is satisfied if the Jacobians $\frac{\partial F_i}{\partial x_i}$ are uniformly bounded.

There exists a communication network among these agents and the network can be described as an undirected or directed graph $\mathcal{G}=(\mathcal{V,E,C})$. Here $\mathcal{V}=\{1,2,\cdots,M\}$ denotes the node set and $\mathcal{E}\subset \mathcal{V}\times \mathcal{V}$ denotes the edge set. $\mathcal{A}=[a_{ij}]\in R^{M\times M}$ is the adjacency matrix. In this framework, if there exists  a connection between $i$ and $j$ nodes(agents), then $a_{ij}>0$; otherwise, $a_{ij}=0$. We assume there is no self-circle in the graph $\mathcal{G}$, i.e., $a_{ii}=0$. A path is a sequence of connected edges in a graph. If there is a path between any two nodes, the graph is said to be connected. If $\mathcal{A}$ is a symmetric matrix,  $\mathcal{G}$ is called an undirected graph. The set of neighbors of node $i$ is denoted by $\mathcal{N}_i=\{j|(i,j)\subset \mathcal{E}\}$.  The in-degree of agent $i$ is denoted as $\text{deg}_i=\sum^{M-1}_{j=1}a_{ij}$ and the degree matrix is denoted as $\mathcal{D}=\text{diag}(\text{deg}_1,\cdots,\text{deg}_M)$. The Laplacian matrix of $\mathcal{G}$ is described as $\mathcal{L=D-A}$.

\textbf{Definition-1:} The nonlinear multi-agent system in \eqref{model} with a given network topology $\mathcal{G}$, and under a distributed control protocol $\mathbf{u}_i=\mu(\mathbf{x}_i,\mathbf{x}_{-i})$, is said to achieve consensus if,
\begin{equation}
\lim_{t\rightarrow \infty}=\|\mathbf{x}_i(t)-\mathbf{x}_j(t)\|=0, \, j=1,\cdots,M,
\end{equation}
where $\mathbf{x}_{-i}$ are the collection of agent $i$'s neighbor's states, i.e., $\mathbf{x}_{-i}=\{\mathbf{x}_{j},j\in \mathcal{N}_i\}$.

With this in mind, in this specific work, we are interested in designing the distributed control strategy $\mathbf{u}_i=\mu(\mathbf{x}_i,\mathbf{x}_{-i})$ using the real-time nonlinear receding horizon control methodology for each agent $i$ to achieve consensus, within the given network topology.

\section{DISTRIBUTED NONLINEAR RECEDING HORIZON CONTROL PROTOCOL}\label{sec:dist_control}

In this presented framework, for each agent $i$, the following optimization problem is utilized to generate the consensus protocol locally, within the network:

\textbf{Problem-1:}
\begin{equation}
\mathbf{u}^*_i(t)=\argmin_{\mathbf{u}_i(t)}J_i(\mathbf{x}_i(t),\mathbf{u}_i(t),\mathbf{x}_{-i}(t))
\end{equation}
subject to
\begin{equation*}
\dot{\mathbf{x}}_i(t)=\mathbf{F}(\mathbf{x}_i(t))+\mathbf{u}_i(t),
\end{equation*}
The performance index is designed as follows:
\begin{equation}\label{cost1}
\begin{split}
J_i&=\varphi_i+\frac{1}{2}\int_t^{t+T}L_i(\mathbf{x}_i,\mathbf{x}_{-i},\mathbf{u}_i),\\
&=\sum_{j\in \mathcal{N}_i} a_{ij}\|\mathbf{x}_i(t+T)-\mathbf{x}_j(t+T)\|^2_{Q_{iN}}+\frac{1}{2}\int_t^{t+T}(\sum_{j\in \mathcal{N}_i} a_{ij}\|\mathbf{x}_i(\tau)-\mathbf{x}_j(\tau)\|^2_{Q_{i}}+\|\mathbf{u}_i\|^2_{R_i}){\rm d}\tau,
\end{split}
\end{equation}
where $Q_{iN}>0$,  $Q_{i}>0$ and $R_{i}>0$ are symmetric matrices, and $T$ is the horizon. We denote $\varphi_i$ to describe the terminal cost for each agent.

Here, we utilize a control scheme to be able to deal with the nonlinear nature of the graph under scrutiny. With this, it is desired to solve the nonlinear optimization problem directly, in real-time.

With the construction of the cost function, as given in \eqref{cost1}, the consensus problem is converted into an optimization procedure. For this purpose, we utilize the powerful nature of \emph{real-time nonlinear receding horizon control} algorithm to generate the distributed consensus protocol by minimizing the associated cost function. In this context, each agent only needs to obtain its neighbors' information once via the given network which is more efficient than the centralized control strategy (and the other distributed strategies) that involve multiple information exchanges and predicted trajectories of states. The performance index evaluates the performance from the present time $t$ to the finite future $t+T$, and then is minimized for each time segment $t$ starting from $\mathbf{x}_i(t)$. With this structure, it is possible to convert the present receding horizon control problem into a family of finite horizon optimal control problems on the artificial $\tau$ axis parameterized by time $t$.

According to the first-order necessary conditions of optimality (i.e. for $\delta J_i=0$), a \emph{local} two-point boundary-value problem (TPBVP)  \cite{bryson1975} is formed as follows:
\begin{equation}\label{tpbvp}
\begin{split}
&{\mathbf{\Lambda}_i}^*_{\tau}(\tau,t)=-H^T_{\mathbf{x}_i},\\
&\mathbf{\Lambda}_i^*(T,t)=\varphi_{\mathbf{x}_i}^T[\mathbf{x}^*(T,t)], \\
&{\mathbf{x}_i}^*_{\tau}(\tau,t)=H^T_{\mathbf{\Lambda}_i},\mathbf{x}_i^*(0,t)=\mathbf{x}_i(t),\\
&H_{\mathbf{u}_i}=0.
\end{split}
\end{equation}
where ${\mathbf{\Lambda}_i}$ denotes the costate of each agent $i$ and $H_i$ is the Hamiltonian which is defined as
\begin{equation}\label{eq:tpbvp_ham}
\begin{split}
H_i&= L_i + \mathbf{\Lambda}_{i}^{*T}\dot{\mathbf{x}}_i\\
&=\frac{1}{2}(\sum_{j\in \mathcal{N}_i} a_{ij}\|\mathbf{x}_i(\tau)-\mathbf{x}_j(\tau)\|^2_{Q_{i}}+\|\mathbf{u}_i\|^2_{R_i})+ \mathbf{\Lambda}_{i}^{*T}(\mathbf{F}(x_i)+\mathbf{u}_i).
\end{split}
\end{equation}
Then we have
\begin{equation}\label{tpbvp_l}
\begin{split}
&{\mathbf{\Lambda}_i}^*_{\tau}(\tau,t)=-[\sum_{j\in \mathcal{N}_i} Q_{i}a_{ij}(\mathbf{x}_i(\tau)-\mathbf{x}_j(\tau))+\mathbf{\Lambda}^T_i\mathbf{F}_{\mathbf{x}_i}(\mathbf{x}_i)],\\
&\mathbf{\Lambda}_i^*(T,t)=\sum_{j\in \mathcal{N}_i} Q_{iN}a_{ij}(\mathbf{x}_i(\tau)-\mathbf{x}_j(\tau)).
\end{split}
\end{equation}
In \eqref{tpbvp}-\eqref{tpbvp_l}, $(~~)^*$ denotes a variable in the optimal control problem so as to distinguish it from its correspondence in the original problem. In this notation, $H_{\mathbf{x}_i}$ denotes the partial derivative of $H$ with respect to ${\mathbf{x}_i}$, and so on.

In this methodology, since the state and costate at $\tau=T$ are determined by the TPBVP in Eq.\eqref{tpbvp} from the state and costate at $\tau=0$, the TPBVP can be regarded as a nonlinear  algebraic equation with respect to the costate at $\tau=0$ as
\begin{equation}\label{p}
\mathbf{P}_i(\mathbf{\Lambda}_i(t),\mathbf{x}_i(t),T,t)=\mathbf{\Lambda}_i^*(T,t)-\varphi_{\mathbf{x}_i}^T[\mathbf{x}_i^*(T,t)]=0,
\end{equation}
where $\mathbf{\Lambda}_i(t)$ denotes the costate at $\tau=0$. The actual local control input for each agent is then given by
\begin{equation}\label{u}
\mathbf{u}_i(t)=\text{arg}\{H_{\mathbf{u}_i}[\mathbf{x}_i(t),\mathbf{\Lambda}_i(t),\mathbf{u}_i(t)]=0\}.
\end{equation}

In this formulation, the optimal control $\mathbf{u}_i(t)$ can be calculated from Eq.\eqref{u} based on $\mathbf{x}_i(t)$ and $\mathbf{\Lambda}_i(t)$ where the ordinary differential equation of $\lambda(t)$ can be solved numerically from Eq.\eqref{p}, in real-time, without any need of an iterative optimization routine. Since the nonlinear equation $\mathbf{P}_i(\mathbf{\Lambda}_i(t),\mathbf{x}_i(t),T,t)$ has to be satisfied at any time $t$, $\frac{{\rm d}\mathbf{P}_i}{{\rm d}t}=0$ holds along the trajectory of the closed-loop system of the receding horizon control. If $T$ is a smooth function of time $t$, the solution of $\mathbf{P}_i(\mathbf{\Lambda}_i(t),\mathbf{x}_i(t),T,t)$ can be tracked with respect to time. However, numerical errors associated with the solution may accumulate as the integration proceeds in practice, and therefore some correction techniques are required to correct such errors in the solution. To address this problem, a stabilized continuation method \cite{kabamba1987,ohtsuka1994,ohtsuka1997,ohtsuka1998} is used. According to this method, it is possible to rewrite the statement as
\begin{equation}\label{cm}
\frac{{\rm d}\mathbf{P}_i}{{\rm d}t}=A_s\mathbf{P}_i,
\end{equation}
where $A_s$ denotes a stable matrix to make the solution converge to zero exponentially.

To evaluate the optimal control by computing derivative of $\mathbf{\Lambda}_i(t)=\mathbf{\Lambda}_i^*(0,t)$ in real time, we consider the partial differentation of \eqref{tpbvp} with respect to time $t$,

\begin{equation}
\begin{split}
&\delta \dot{\mathbf{x}}_i=f_{\mathbf{x}_i}\delta \mathbf{x}_i+f_{\mathbf{u}_i}\delta \mathbf{u}_i,\\
&\delta \dot{\mathbf{\Lambda}}_i=-H_{\mathbf{x}_i\mathbf{x}_i}\delta \mathbf{x}_i-H_{\mathbf{x}_i\mathbf{\Lambda}_i}\delta \mathbf{\Lambda}_i-H_{\mathbf{x}_i\mathbf{u}_i}\delta \mathbf{u}_i\\
&0=H_{\mathbf{u}_i\mathbf{x}_i}\delta \mathbf{x}_i+f_{\mathbf{u}_i}^T\delta \mathbf{\Lambda}_i+H_{\mathbf{u}_i\mathbf{u}_i}\delta \mathbf{u}_i.
\end{split}
\end{equation}

Since $\delta \mathbf{u}_i=-H^{-1}_{\mathbf{u}_i\mathbf{u}_i}(H_{\mathbf{u}_i\mathbf{x}_i}\delta \mathbf{x}_i+f_{\mathbf{u}_i}^T\delta \mathbf{\Lambda}_i)$, we have
\begin{align*}
&\delta \dot{\mathbf{x}}_i=(f_{\mathbf{x}_i}-f_{\mathbf{u}_i}H^{-1}_{\mathbf{u}_i\mathbf{u}_i}H_{\mathbf{u}_i\mathbf{u}_i})\delta \mathbf{x}_i-f_{\mathbf{u}_i}H^{-1}_{\mathbf{u}_i\mathbf{u}_i}f_{\mathbf{u}_i}^T\delta \mathbf{\Lambda}_i,\\
&\delta \dot{\mathbf{\Lambda}}_i=-(H_{\mathbf{x}_i\mathbf{x}_i}-H_{\mathbf{x}_i\mathbf{u}_i}H^{-1}_{\mathbf{u}_i\mathbf{u}_i}H_{\mathbf{u}_i\mathbf{x}_i})\delta \mathbf{x}_i-(f_{\mathbf{x}_i}^T-H_{\mathbf{x}_i\mathbf{u}_i}H^{-1}_{\mathbf{u}_i\mathbf{u}_i}f_{\mathbf{u}_i}^T)\delta \mathbf{\Lambda}_i,
\end{align*}  

which leads to the following form of a linear differential equation:
\begin{equation}\label{pd}
\frac {\partial}{\partial \tau}\begin{bmatrix}{\mathbf{x}_i}^*_t-{\mathbf{x}_i}^*_{\tau}\\
{\mathbf{\Lambda}_i}^*_t-{\mathbf{\Lambda}_i}^*_{\tau} \end{bmatrix}=\begin{bmatrix}A_i&-B_i\\
-C_i&-A_i^T \end{bmatrix}\begin{bmatrix}{\mathbf{x}_i}^*_t-{\mathbf{x}_i}^*_{\tau}\\
{\mathbf{\Lambda}_i}^*_t-{\mathbf{\Lambda}_i}^*_{\tau} \end{bmatrix}
\end{equation}
where
$$A_i=f_{\mathbf{x}_i}-f_{\mathbf{u}_i}H^{-1}_{{\mathbf{u}_i}{\mathbf{u}_i}}H_{{\mathbf{u}_i}{\mathbf{x}_i}},$$
$$B_i=f_{\mathbf{u}_i}H^{-1}_{{\mathbf{u}_i}{\mathbf{u}_i}}f_{{\mathbf{u}_i}}^T,$$
$$C_i=H_{{\mathbf{x}_i}{\mathbf{x}_i}}-H_{{\mathbf{x}_i}{\mathbf{u}_i}}H^{-1}_{{\mathbf{u}_i}{\mathbf{u}_i}}H_{{\mathbf{u}_i}{\mathbf{x}_i}}.$$ And the matrix $H_{\mathbf{u}_i\mathbf{u}_i}$ should be nonsingular.

In order to reduce the computational cost without resorting to any approximation technique, the backward-sweep method \cite{bryson1975,ohtsuka1997,ohtsuka1998} is implemented. The derivative of the function $\mathbf{P}_i$ with respect to time is given by
\begin{equation}\label{df}
\begin{split}
\frac{{\rm d}\mathbf{P}_i}{{\rm d}t}=&{\mathbf{\Lambda}_i}^*_t(T,t)-\varphi_{{\mathbf{x}_i}{\mathbf{x}_i}}{\mathbf{x}_i}^*_t(T,t)+[{\mathbf{\Lambda}_i}^*_{\tau}(T,t)-\varphi_{{\mathbf{x}_i}{\mathbf{x}_i}}{\mathbf{x}_i}^*_{\tau}(T,t)]\frac{{\rm d}T}{{\rm d}t},
\end{split}
\end{equation}
where $x^*_{\tau}$ and ${\mathbf{\Lambda}_i}^*_{\tau}$ are given by \eqref{tpbvp}.

The relationship between the costate and other variables is assumed as follows:
\begin{equation}\label{relation}
{\mathbf{\Lambda}_i}^*_t-{\mathbf{\Lambda}_i}^*_\tau=S_i(\tau,t)({\mathbf{x}_i}^*_t-{\mathbf{x}_i}^*_\tau)+c_i(\tau,t).
\end{equation}
So we have
\begin{equation}\label{sct}
\begin{split}
&S_i(T,t)=\varphi_{\mathbf{x}_i\mathbf{x}_i}\mid_{\tau=T},\\
&c_i(T,t)=(H_{\mathbf{x}_i}^T+\varphi_{\mathbf{x}_i\mathbf{x}_i}f)\mid_{\tau=T}(1+\frac{{\rm d}T}{{\rm d}t})+A_s\mathbf{P}_i.
\end{split}
\end{equation}
Then according to \eqref{relation} and \eqref{pd}, we have the following differential equations:
\begin{equation}\label{sc}
\begin{split}
&\frac {\partial {S_i}}{\partial \tau}=-A_i^TS_i-S_iA_i+S_iB_iS_i-C_i,\\
&\frac {\partial {c_i}}{\partial \tau}=-(A_i^T-S_iB_i)c_i.
\end{split}
\end{equation}
Based on \eqref{relation}, the differential equation of the costate to be integrated in real time is obtained as:
\begin{align}\label{dl}
\frac{{\rm d}\mathbf{\Lambda}_i(t)}{{\rm d}t}=-H_{\mathbf{x}_i}^T+c_i(0,t).
\end{align}

Here, at each time $t$, the TPBVP equations are integrated forward along the $\tau$ axis, and then \eqref{sc} are integrated backward with terminal conditions provided in \eqref{sct}. Next, the differential equation of $\mathbf{\Lambda}_i(t)$ is integrated for one step along the $t$ axis so as to update the local optimal controller for each agent. If the matrix $H_{\mathbf{u}_i\mathbf{u}_i}$ is nonsingular, the algorithm is executable regardless of controllability or stabilizability of the system.

\textbf{Lemma-1:} The cost function, defined in Eq. \eqref{cost1}, is strictly convex and guarantees the global minimum.

\emph{Proof}: Since all weighting functions maintain positive definite nature in their structure, such as $Q_{iN}>0$,  $Q_{i}>0$, $R_{i}>0$, from the Karush-Kuhn-Tucker(KKT) conditions \cite{antoniou2007}, the proposed method gurantees the global minima. \qedsymbol

\section{CONVERGENCE AND STABILITY ANALYSIS}\label{sec:stab_analysis}
For the sake of clarity, and without loss of generality, here we define the consensus error as
$$\delta_1(t)=\mathbf{x}_i(t)-\mathbf{x}_1(t)$$
for all $i$ and $\Delta(t)=\text{col}(\delta_1(t),\delta_2(t),\cdots,\delta_M(t))$. The optimal control $\mathbf{U}$ is denoted as $\mathbf{U}=\text{col}(\mathbf{u}_1,\cdots,\mathbf{u}_M)$ for all $i$. 

The cost function can be written as
\begin{equation}
J=\mathbf{\Phi}+\frac{1}{2}\int_t^{t+T}[\Delta^{*T}Q\Delta^*+\mathbf{U}^{*T}R\mathbf{U}^*]{\rm d}\tau,
\end{equation}
where $\mathbf{\Phi}=\sum_{i \in M}\varphi_i$, $Q=\text{col}(Q_1, \cdots, Q_M)$ and $R=\text{col}(R_1, \cdots, R_M)$.

In order to ensure the closed-loop stability of the proposed nonlinear receding horizon control scheme, we first consider the case that terminal cost $\mathbf{\Phi}=0$ and introduce following definitions.

In this regard, we assume the sublevel sets 
\begin{equation*}
\Gamma_r^{\infty}=\{\Delta\in\Gamma^{\infty}:J^*_{\infty}<r^2\}
\end{equation*}
are compact and path connected where $J^*_{\infty}=\int_0^{\infty}[\Delta^{*T}Q\Delta^*+\mathbf{U}^{*T}R\mathbf{U}^*]{\rm d}\tau$ and moreover $\Gamma^{\infty}=\cup_{r\ge 0}\Gamma_r^{\infty}$. We use $r^2$ here to reflect the fact that the cost function is quadratically bounded. And therefore the sublevel set of $\Gamma_r^{T}=\{\Delta\in\Gamma^{\infty}:J^*_{T}<r^2\}$ where $J^*_{T}=\int_t^{t+T}[\Delta^{*T}Q\Delta^*+\mathbf{U}^{*T}R\mathbf{U}^*]{\rm d}\tau$.

\textbf{Lemma-2:} (Dini \cite{jadbabaie2005} ) Let $\{f_n\}$ be a sequence of upper semi-continuous, real-valued functions on a countably compact space $X$, and suppose that for each $x\in X$, the sequence $\{f_n(x)\}$ decreases monotonically to zero. Then the convergence is uniform.

\textbf{Theorem-1:} \cite{jadbabaie2005} Let $r$ be given as $r>0$ and suppose that the terminal cost is equal to zero. For each sampling time $t_s >0$, there exists a horizon window  $T^*<\infty$ such that, for any $T \geqslant T^*$, the receding horizon scheme is asymptotically stabilizing. 

\begin{proof}
By the principle of optimality, we have
\begin{equation*}
J^*_T(\Delta)=\int_{t}^{t+t_s}(\Delta_T^{*T}Q\Delta_T^*+\mathbf{U}_T^{*T} R\mathbf{U}_T^*)d\tau+J^*_{T-t_s}(\Delta_T^*)
\end{equation*}
where $t_s \in [t,t+T]$ is the sampling time and $J^*_{t_s}(\Delta)=\int_{t}^{t+t_s}(\Delta_T^{*T}Q\Delta_T^*+\mathbf{U}_T^{*T} R\mathbf{U}_T^*)d\tau$, so that
\begin{equation*}
\begin{split}
J^*_{T-t_s}(\Delta_T^*)-J^*_{T-t_s}(\Delta)&=J^*_{T}(\Delta)-J^*_{T-t_s}(\Delta)-\int_t^{t+t_s}(\Delta_T^{*T}Q\Delta_T^{*}+\hat{\Theta}_T^{*T} R\hat{\Theta}_T^{*})d\tau\\
&\leq J^*_{t_s}(\Delta)+J^*_{T}(\Delta)-J^*_{T-t_s}(\Delta)
\end{split}
\end{equation*}
Since the terminal cost is equal to zero, it is clear that $T_1<T_2$. This implies that $J^*_{T_1}(\Delta)<J^*_{T_2}(\Delta)$ holds for all $\Delta$ so that
\begin{equation*}
J^*_{T-t_s}(\Delta_T^*)-J^*_{T-t_s}(\Delta) \leq J^*_{t_s}(\Delta)+J^*_{\infty}(\Delta)-J^*_{T-t_s}(\Delta).
\end{equation*}
is satisfied.
If we can show, for example, that there exists a $T^*$ such that $T>T^*$ yields into
\begin{equation*}
J^*_{\infty}(\Delta)-J^*_{T-t_s}(\Delta) \leq \frac{1}{2}J^*_{t_s}(\Delta)
\end{equation*}
for all $\Delta \in \Gamma_r^{\infty}$, stability over any sublevel set of $J^*_{T-t_s}(\cdot)$ that is contained in $\Gamma_r^{\infty}$ will be assured. To that end, define, for $\Delta \in \Gamma_r^{\infty}$

\begin{equation*}
\psi_T(\Delta)=
\begin{cases}
\frac{J^*_{\infty}(\Delta)-J^*_{T-t_s}(\Delta)}{J^*_{t_s}(\Delta)}, & \quad \Delta \neq 0 \\
\lim \sup_{x\rightarrow 0} \psi_T(\Delta), & \quad \Delta=0
\end{cases}
\end{equation*}
where $\psi_T(\cdot)$ is upper semicontinuous on $\Gamma_r^{\infty}$. It is clear that $\psi_T(\cdot)$ is a monotonically decreasing family of upper semicontinuous functions defined over the compact set $\Gamma_r^{\infty}$. Thus, by Dini's theorem (as stated in \cite{jadbabaie2005}), there exists a $T^*<\infty$ such that $\psi_T(\Delta)<\frac{1}{2}$ for all $\Delta \in \Gamma_r^{\infty}$ and all $T\ge T^*$. Here, for each $r_1>0$ we have $\Gamma_{r_1}^{T-t_s}\subset \Gamma_r^{\infty}$ satisfied, leading to

\begin{equation*}
J^*_{T-t_s}(\Delta_T^*)-J^*_{T-t_s}(\Delta) \leq -\frac{1}{2}J^*_{t_s}(\Delta)
\end{equation*}
for all $\Delta \in \Gamma_{r_1}^{T-t_s}$.

\end{proof}

Next, we present the closed-loop stability of the proposed nonlinear receding horizon control scheme with locally quadratic terminal cost, i.e. $\mathbf{\Phi}=\sum_{i \in M}\varphi_i$.

\textbf{Theorem-2:} \cite{jadbabaie2005} Let $r$ be given as $r>0$ and suppose that the terminal cost is nonnegative and locally quadratically bounded. For each sampling time $t_s>0$, there exists a horizon window  $T^*<\infty$ such that, for any $T \geqslant T^*$, the receding horizon scheme is asymptotically stabilizing. 

\begin{proof}
We use $J^*_{T,0}(\cdot)$ to denote the cost function with zero terminal cost and $J^*_{T,1}(\cdot)$ to denote the cost function with locally quadratic terminal cost. It is clear to show that
\begin{equation*}
J^*_{T,0}(\Delta)\leq J^*_{T}(\Delta)\leq J^*_{T,1}(\Delta),
\end{equation*} 
and then 
\begin{equation*}
|J^*_{T}(\Delta)- J^*_{\infty}(\Delta)|\leq \max \{J^*_{\infty}(\Delta)-J^*_{T,0}(\Delta),J^*_{T,1}(\Delta)-J^*_{\infty}(\Delta)\},
\end{equation*} 
for all $\Delta \in \Gamma_r^{\infty}$ so that $J^*_{T}(\cdot)$ also converge uniformly to $J^*_{\infty}(\cdot)$ with respect to any locally quadratic positive definite terminal cost. The theorem uses the results of \textbf{Theorem-1}.

\end{proof}

\textbf{Corollary-1:} Consider the nonlinear multi-agent system given in Eq.\eqref{model} and assume that the graph $\mathcal{G}$, defining the topology of the multi agent system, is connected. For the given distributed control protocol in \textbf{Problem-1}, based on \textbf{Theorem-2}, there exists a large enough value of horizon $T$ which guarantees the consensus error $\Delta$ to remain asymptotically stable to achieve consensus in the multi-agent system.

Although the back-ward sweep algorithm is executable whenever the system is stable or not, with this result, when the optimization horizon is chosen to be sufficiently long, the non-increasing monotonicity of the cost function becomes a sufficient condition for the stability. Therefore, we can also ensure the stability of the multi-agent nonlinear system by distributed nonlinear receding horizon control method.

\section{EXTENSION TO LEADER-FOLLOWING CASE}\label{sec:leader-following}

In this section, we extend the proposed scheme to the leader-following consensus problem(s). Assume a multi-agent system consisted of $M$ nonlinear agents and a leader. The dynamics of each agent is given in Eq. \eqref{model} and the leader, indexed as $i=0$, has the following dynamical system structure:
 \begin{equation}\label{model_l}
\dot{\mathbf{x}}_0=f(\mathbf{x}_0,t)=\mathbf{F}(x_0(t)),
\end{equation}
where $\mathbf{x}_{0}=(x_{01},x_{02},\cdots,x_{0n})^{T}$ is the state vector.

To describe the connection between the leader and the followers in $\mathcal{G}$, a leader adjacency matrix $\mathcal{A}_0=\text{diag}(a_{10},\cdots,a_{M0})$, where $a_{i0}>0$ if follower $i$ is connected to the leader, otherwise $a_{i0}=0$. Then a new augmented matrix $\mathcal{H}=\mathcal{A}+\mathcal{A}_0$ is defined.

\textbf{Definition-2:} The nonlinear multi-agent system in \eqref{model} and \eqref{model_l} with a given network topology $\mathcal{G}$, and under a distributed control protocol $\mathbf{u}_i=\mu(\mathbf{x}_0,\mathbf{x}_i,\mathbf{x}_{-i})$, is said to achieve the leader-following consensus such that all $M$ follower agents converge to the leader, if,
\begin{equation}
\lim_{t\rightarrow \infty}=\|\mathbf{x}_i(t)-\mathbf{x}_0(t)\|=0, \, i=1,\cdots,M,
\end{equation}
where $\mathbf{x}_{-i}$ are also the collection of agent $i$'s neighbor's states, i.e., $\mathbf{x}_{-i}=\{\mathbf{x}_{j},j\in \mathcal{N}_i\}$.

According to the presented framework, next, we are interested in designing the distributed control strategy $\mathbf{u}_i=\mu(\mathbf{x}_0, \mathbf{x}_i,\mathbf{x}_{-i})$ using the abovementioned real-time nonlinear receding horizon control methodology for each agent $i$ to achieve leader-following consensus, within the given network topology.

Similarly, for each following agent $i$, the following optimization problem is utilized to generate the consensus protocol locally, within the network:

\textbf{Problem-2:}
\begin{equation}
\mathbf{u}^*_i(t)=\text{argmin}J_i(\mathbf{x}_0, \mathbf{x}_i,\mathbf{u}_i,\mathbf{x}_{-i})
\end{equation}
subject to
\begin{equation*}
\begin{cases}
\dot{\mathbf{x}}_0(t)=\mathbf{F}(x_0(t)),\\
\dot{\mathbf{x}}_i(t)=\mathbf{F}(\mathbf{x}_i(t))+\mathbf{u}_i(t),
\end{cases}
\end{equation*}
The corresponding performance index is designed as follows:
\begin{equation}\label{cost}
\begin{split}
J_i=&\varphi_i+\frac{1}{2}\int_t^{t+T}L_i(\mathbf{x}_i,\mathbf{x}_0,\mathbf{x}_{-i},\mathbf{u}_i),\\
=&\sum_{i\in M} c_{i0}\|\mathbf{x}_i(t+T)-\mathbf{x}_0(t+T)\|^2_{Q_{i0}}+\sum_{j\in \mathcal{N}_i} c_{ij}\|\mathbf{x}_i(t+T)-\mathbf{x}_j(t+T)\|^2_{Q_{iN}}\\
&+\frac{1}{2}\int_t^{t+T}(\sum_{j\in \mathcal{N}_i} c_{ij}\|\mathbf{x}_i(\tau)-\mathbf{x}_j(\tau)\|^2_{Q_{i}}+\|\mathbf{u}_i\|^2_{R_i}){\rm d}\tau,
\end{split}
\end{equation}
where the terminal cost $\varphi_i=\sum_{i\in M} c_{i0}\|\mathbf{x}_i(t+T)-\mathbf{x}_0(t+T)\|^2_{Q_{i0}}+\sum_{j\in \mathcal{N}_i} c_{ij}\|\mathbf{x}_i(t+T)-\mathbf{x}_j(t+T)\|^2_{Q_{iN}}$, $Q_{i0}>0$, $Q_{iN}>0$,  $Q_{i}>0$ and $R_{i}>0$ are symmetric matrices, and $T$ is also the horizon.

We utilize the same framework and algorithm as in Section-\ref{sec:dist_control} to solve the nonlinear leader-following consensus problem directly. In addition, when the optimization horizon is chosen to be sufficiently long, the non-increasing monotonicity of the cost function becomes a sufficient condition for the stability.

\section{NUMERICAL EXAMPLE AND SIMULATION RESULTS}\label{sec:simulation}
In this section, we demonstrate the validity and feasibility of proposed scheme on several multi-agent nonlinear \emph{chaotic} systems.

\begin{figure*}
 \centering
    \includegraphics[width=20cm]{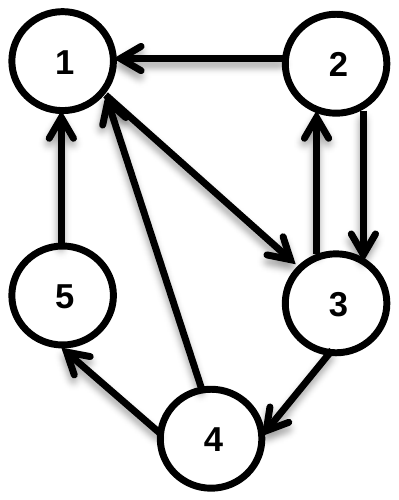}
    \vspace*{-20.1cm}
    \caption{ \label{t1} The communication topology of the multi-agent Lorenz chaotic systems with $5$ agents.}
 \end{figure*}

\textbf{Example-1:} First, consider a multi-agent system with $5$ agents (as shown in Fig. \ref{t1}), where each agent is modeled as the Lorenz \emph{chaotic} system \cite{lorenz1963}:
\begin{equation}\label{lorenz}
\begin{cases}
\dot{x}_{i1}=10(x_{i2}-x_{i1}), \\
\dot{x}_{i2}=28x_{i1}-x_{i1}x_{i3}-x_{i2},\\
\dot{x}_{i3}=x_{i1}x_{i2}-\frac{8}{3}x_{i3},
\end{cases}
\end{equation}
where $x_i=(x_{i1},x_{i2},x_{i3})^{T}$ are the state of the $i$-th agent. The adjacency matrix $\mathcal{A}$ of $\mathcal{G}$ is given as
\begin{equation}
\begin{bmatrix}0&0&1&0&0 \\ 1&0&1&0&0\\0&1&0&1&0\\1&0&0&0&1\\1&0&0&0&0\end{bmatrix}
\end{equation}

Here, the initial states are given by
\begin{equation}\label{initial1}
\begin{bmatrix}x_{i1}(0) \\ x_{i2}(0)\\x_{i3}(0)\end{bmatrix}=\begin{bmatrix}1&2&-10&9\\ 10& -1& 20& -10\\2& 5& 8& -2\end{bmatrix},
\end{equation}

It can be seen in Fig. \ref{t1} that the directed graph $\mathcal{G}$  of Lorenz network is connected. The weighting matrices in the cost function are designed as $Q_{iN}=Q_{i}=R_{i}=\text{diag}(1,1,1)$ for all agents. The stable matrix is designed as $A_s=-50I$.

The horizon $T$ in the performance index is given by
\begin{equation}\label{h_t}
T(t)=T_f(1-e^{-\alpha t}),
\end{equation}
where $T_f=1$ and $\alpha=0.01$.

\begin{figure}
 \centering
    \includegraphics[width=15cm]{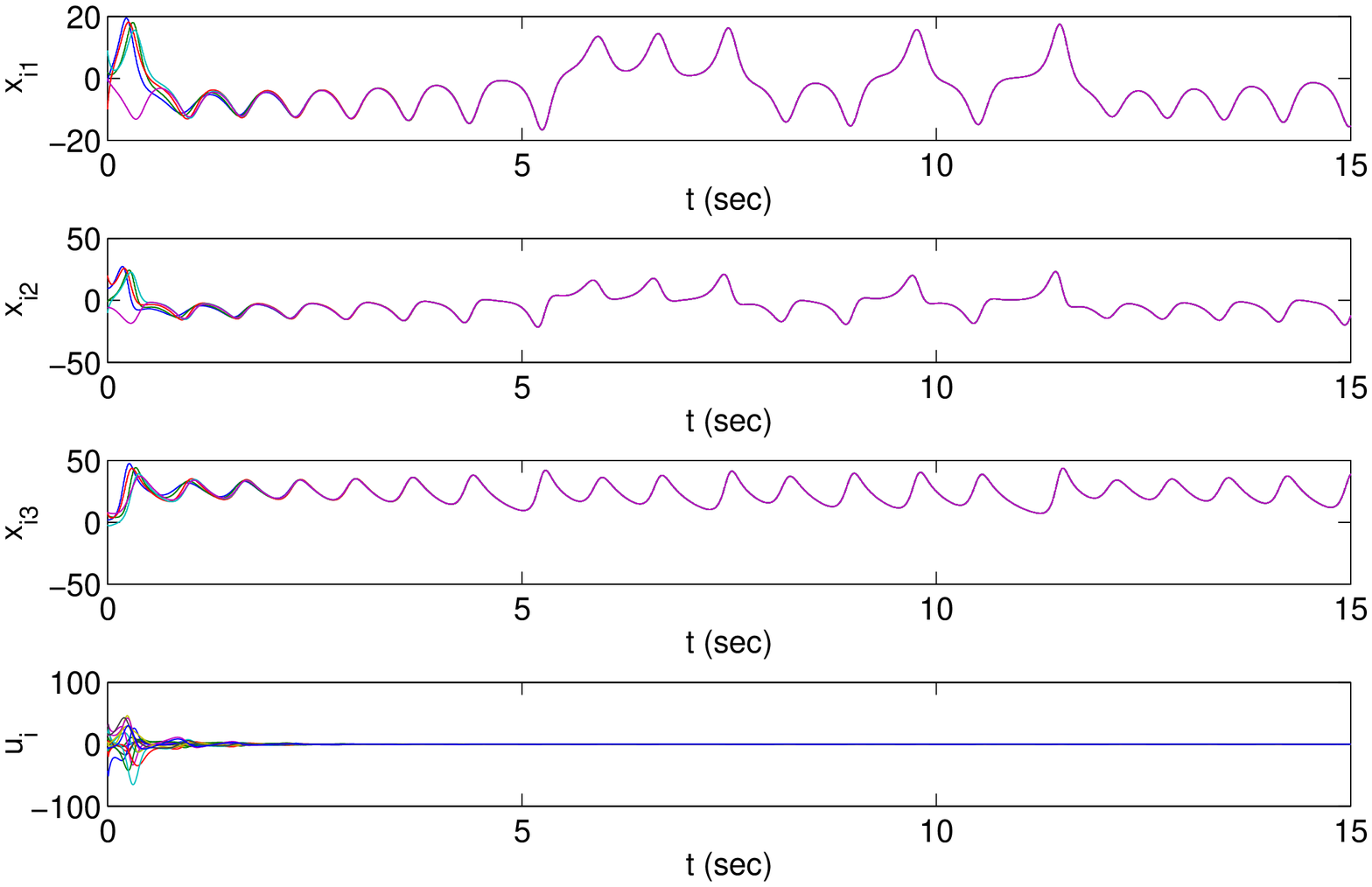}
    \vspace*{-0.5cm}
    \caption{\label{fig_m1} The trajectories of all agents $x_i(t) (i=1,\cdots,5)$ of Lorenz chaotic system and all control protocol $u_i$ generated by distributed NRHC.}
 \end{figure}

The simulation is implemented in MATLAB, where the time step on the $t$ axis
is $0.01s$ and the time step on the artificial $\tau$ axis is $0.005s$. Fig. \ref{fig_m1} depicts the trajectories of this multi-agent Lorenz systems with initial conditions defined in Eq. \eqref{initial1}, where it is possible to observe that the proposed distributed real-time nonlinear receding horizon control strategy results in a consensus, clearly demonstrating the effectiveness of the algorithm. Here, the horizon length is kept sufficiently long to ensure the stability.

\begin{figure}
 \centering
    \includegraphics[width=20cm]{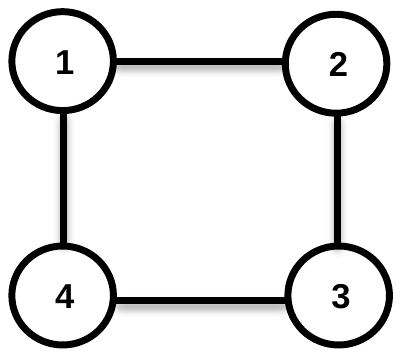}
    \vspace*{-21.1cm}
    \caption{ \label{t2} The communication topology of the multi-agent L\"{u} chaotic systems with $4$ agents.}
 \end{figure}

\textbf{Example-2:} Next, consider a multi-agent system with $4$ agents (as shown in Fig. \ref{t2}), where each agent is modeled as the L\"{u} chaotic system \cite{lu_chen2002}:
\begin{equation}\label{lv}
\begin{cases}
\dot{x}_{i1}=36(x_{i2}-x_{i1}), \\
\dot{x}_{i2}=-x_{i1}x_{i3}+13x_{i2},\\
\dot{x}_{i3}=x_{i1}x_{i2}-3x_{i3},
\end{cases}
\end{equation}
Here, $x_i=(x_{i1},x_{i2},x_{i3})^{T}$ represents the state of the $i$-th agent. In this case, the adjacency matrix $\mathcal{A}$ of $\mathcal{G}$ is given as
\begin{equation}
\begin{bmatrix}0&1&0&1 \\ 1&0&1&0\\0&1&0&1\\1&0&1&0\end{bmatrix}
\end{equation}
It can be seen in Fig. \ref{t2} that the undirected graph $\mathcal{G}$  of L\"{u} network is connected.  The weighting matrices in the cost function are designed as $Q_{iN}=Q_{i}=R_{i}=\text{diag}(1,1,1)$ for all agents.The stable matrix is designed as $A_s=-50I$. The horizon $T$ in the performance index is given in \eqref{h_t} with $T_f=1$ and $\alpha=0.01$.

\begin{figure}
 \centering
    \includegraphics[width=15cm]{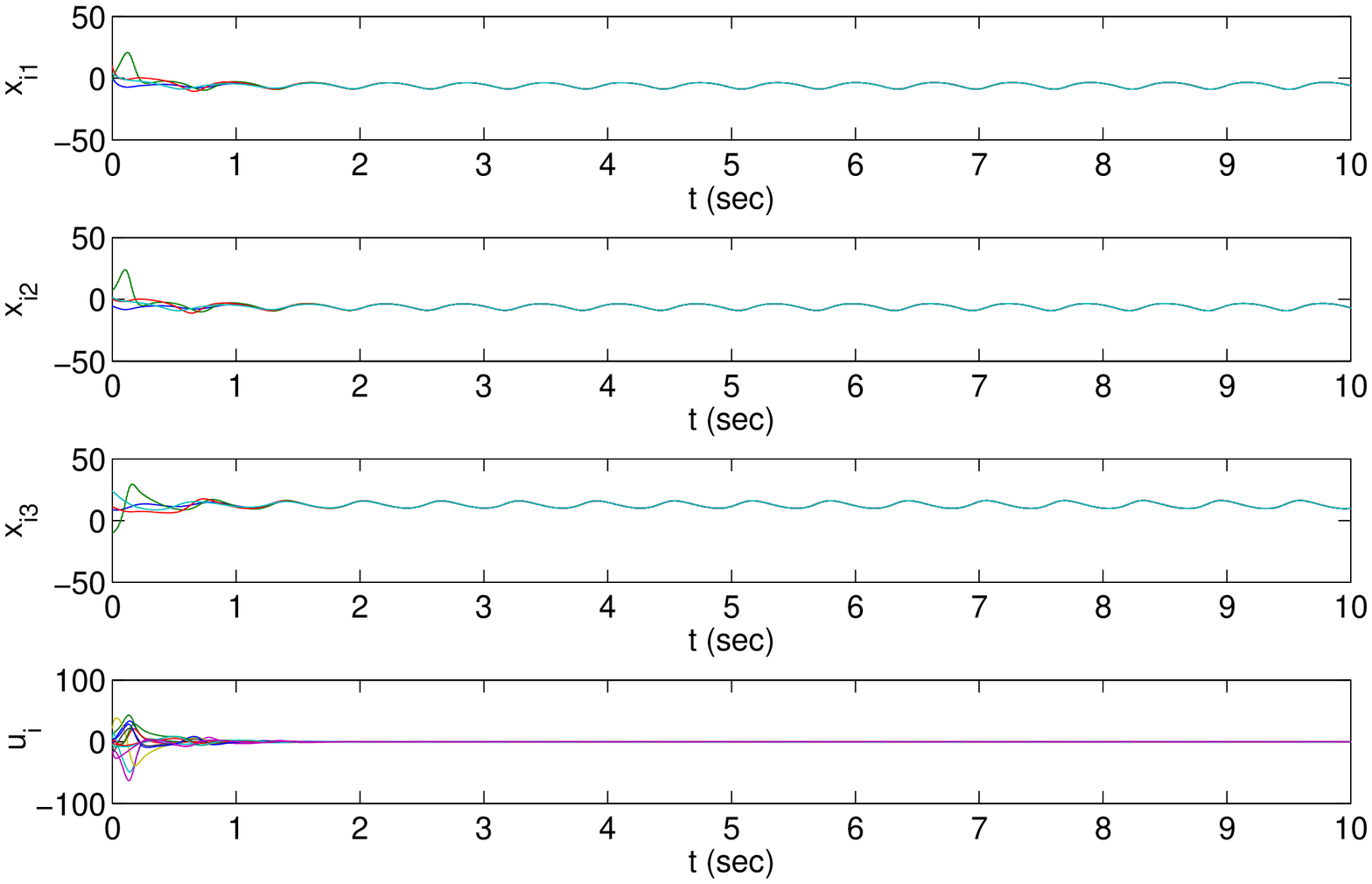}
    \vspace*{-0.5cm}
    \caption{ \label{fig_m2} The trajectories of all agents $x_i(t) (i=1,\cdots,4)$ of L\"{u} chaotic system and all control protocol $u_i$ generated by distributed NRHC.}
 \end{figure}

The simulation is implemented in MATLAB, where the time step on the $t$ axis
is $0.01s$ and the time step on the artificial $\tau$ axis is $0.005s$. Fig. \ref{fig_m2} depicts the trajectories of this multi-agent L\"{u} systems in \eqref{lv} with initial conditions defined in Eq. \eqref{initial1} and the corresponding control strategies, where the consensus is achieved through the suggested distributed real-time nonlinear receding horizon control method.

\textbf{Example-3:} Consider a leader-following system with one leader and $4$ following agents, where each agent and the leader are also modeled as the Lorenz chaotic system (as shown in Eq. \eqref{lorenz}). The augmented adjacency matrix $\mathcal{H}$ of $\mathcal{G}$ is given as
\begin{equation}
\begin{bmatrix}1&1&0&0 \\ 1&0&1&0\\0&1&1&1\\0&0&1&0\end{bmatrix}
\end{equation}
where the diagonal elements show the connection between the leader and the followers. The weighting matrices in the cost function are designed as $Q_{i0}=\text{diag}(10,10,10)$ and $Q_{iN}=Q_{i}=R_{i}=\text{diag}(1,1,1)$ for all agents. The stable matrix is designed as $A_s=-50I$.

The initial states are given by
\begin{equation}\label{initial2}
\begin{bmatrix}x_{i1}(0) \\ x_{i2}(0)\\x_{i3}(0)\end{bmatrix}=\begin{bmatrix}0.1 &-1&2&-10&9\\ 0.2&10& -1& 20& -10\\0.3&2& 5& 8& -2\end{bmatrix},
\end{equation}
where $i=0,1,\cdots,4$.

The horizon $T$ in the performance index is given by
\begin{equation}\label{horizon}
T(t)=T_f(1-e^{-\alpha t}),
\end{equation}
where $Tf=1$ and $\alpha=0.01$.

\begin{figure}
\centering
 \vspace*{-.1cm}
    \includegraphics[width=15cm]{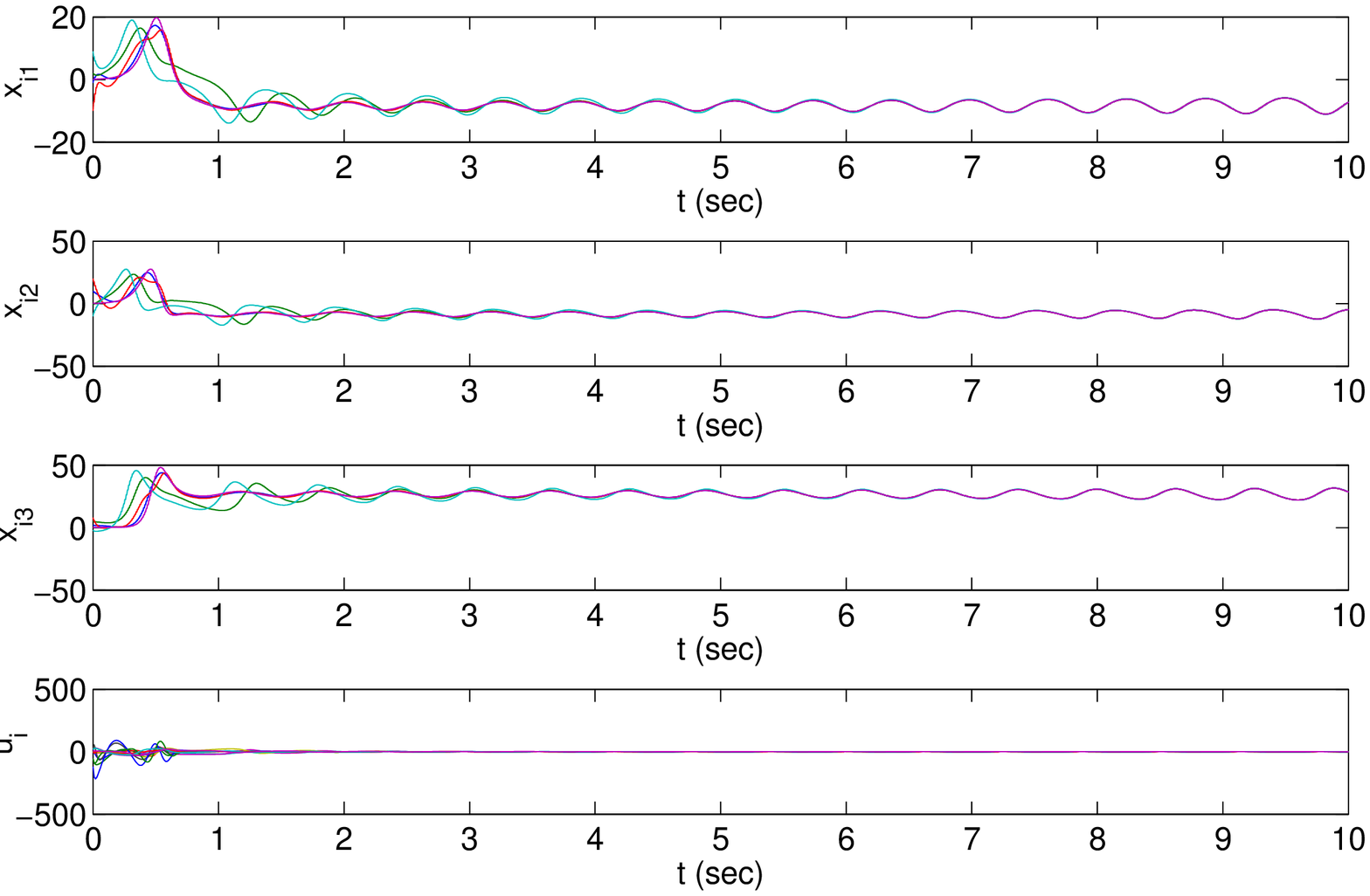}
    \vspace*{-.5cm}
    \caption{\label{fig_lf1} The trajectories of leader and following agents $x_i(t) (i=0,1,\cdots,4)$ of Lorenz chaotic system and all control protocol $u_i$ generated by distributed NRHC.}
 \end{figure}

The simulation is implemented in MATLAB, where the time step on the $t$ axis
is $0.01s$ and the time step on the artificial $\tau$ axis is $0.005s$. Fig. \ref{fig_lf1} depicts the trajectories
of this leader-following Lorenz system with initial condition in Eq. \eqref{initial2} and the corresponding control strategies where the leader-following consensus is reached by using the distributed real-time nonlinear receding horizon control method.

\textbf{Example-4:} We consider the leader-following multi-agent system with one leader agent and $4$ agents and each agent is modeled as the Chen chaotic system \cite{chen1999}:
\begin{equation}\label{chen}
\begin{cases}
\dot{x}_{i1}=35(x_{i2}-x_{i1}), \\
\dot{x}_{i2}=-7x_{i1}-x_{i1}x_{i3}+28x_{i2},\\
\dot{x}_{i3}=x_{i1}x_{i2}-3x_{i3},
\end{cases}
\end{equation}
where $x_i=(x_{i1},x_{i2},x_{i3})^{T}$ are the state of the $i$-th agent. The adjacency matrix $\mathcal{H}$ of $\mathcal{G}$ is given as
\begin{equation}
\begin{bmatrix}1&1&1&1 \\ 1&0&1&0\\0&0&0&1\\1&0&0&1\end{bmatrix}
\end{equation}

\begin{figure}
 \centering
 \vspace*{-1.cm}
    \includegraphics[width=15cm]{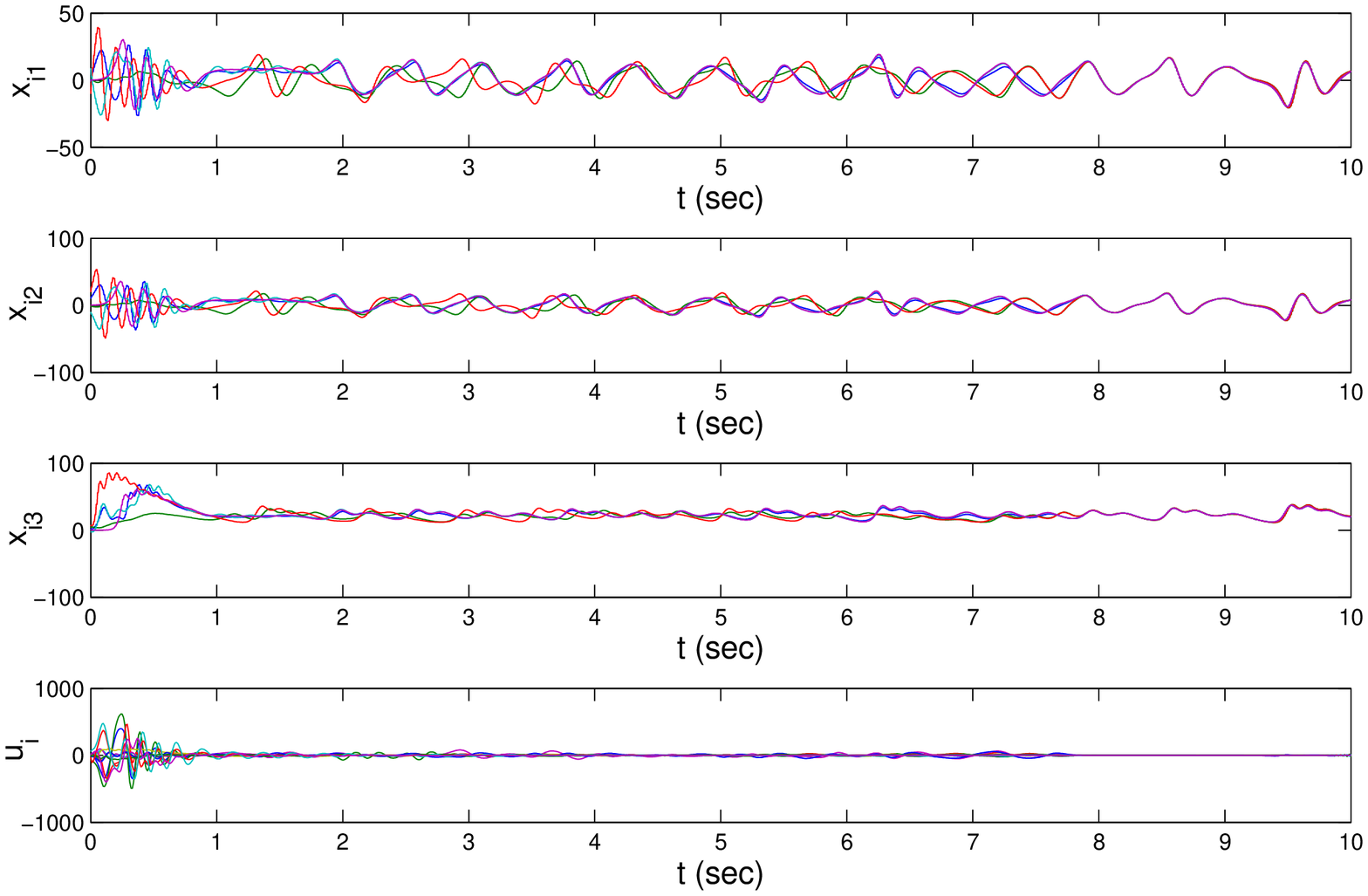}
    \vspace*{-.5cm}
    \caption{\label{fig_lf2}The trajectories of leader and following agents $x_i(t) (i=0,1,\cdots,4)$ of Chen chaotic system and all control protocol $u_i$ generated by distributed NRHC.}
 \end{figure}

The weighting matrices in the cost function are designed as $Q_{i0}=\text{diag}(10,10,10)$ and $Q_{iN}=Q_{i}=R_{i}=\text{diag}(1,1,1)$ for all agents. The stable matrix is designed as $A_s=-50I$. The horizon $T$ in the performance index is given in \eqref{horizon} with $T_f=0.5$ and $\alpha=0.01$. Fig. \ref{fig_lf2} depicts the trajectories
of this leader-following Chen system in \eqref{chen} with initial condition in Eq. \eqref{initial2} and the corresponding control strategies where the consensus is reached by using the distributed real-time nonlinear receding horizon control method.

As it can be seen easily from above example, proposed real-time nonlinear receding horizon control methodology is working remarkable within given network topology and graph. All the agents in the systems are able to reach consensus. Here, again, the horizon length is kept sufficiently long to ensure the stability.

\section{CONCLUSIONS AND FUTURE WORKS}\label{sec:conclusions}

In this paper, we investigated the multi-agent consensus problem of nonlinear systems by using distributed real-time nonlinear receding horizon control methodology. Different from the previous works, we solved the nonlinear optimal consensus problem directly, without any need or the utilization of linearization techniques and/or iterative procedures. Based on the stabilized continuation method, the backward sweep algorithm is implemented to minimize the consensus error among the agents and the local control strategy is integrated in real time. We provided stability guarantees of the systems if the horizon length is kept sufficiently long. Several benchmark examples with different topologies demonstrates the applicability and significant outcomes of proposed scheme on nonlinear chaotic systems. For future studies, it is the authors intention to extend the proposed method to the systems with switching topologies and embedded communication time-delays.

%
%
%
%

\newpage

\begin{thebibliography}{99}

 \bibitem{saber2004} R. Olfati-Saber, R. M. Murray, Consensus problems in networks of agents with switching topology and time-delays,  IEEE Trans. on Automat. Contr. 49 (2004)  1520-1533.

\bibitem{dunbar2006} W. Dunbar, R. M. Murray, Distributed receding horizon control for multi-vehicle formation stabilization, Automatica 42 (2006) 549–558.

\bibitem{dunbar2007} W. Dunbar, Distributed receding horizon control of dynamically coupled nonlinear systems, IEEE Trans. on Automat. Contr. 52 (2007) 1249-1263.

\bibitem{hui2008} Q. Hui, W. Haddad, Distributed nonlinear control algorithms for network consensus, Automatica 44 (2008) 2375-2381.

\bibitem{keviczky2008} T. Keviczky, F. Borrelli, K. Fregene, D. Godbole, G. J. Balas, Decentralized receding horizon control and coordination
of autonomous vehicle formations, IEEE Trans. on Control Systems Technology 16 (2008) 19–33.

\bibitem{you2011} K. You, L. Xie, Network topology and communication data rate for consensusability of discrete-time multi-agent systems,  IEEE Trans. on Automat. Contr. 56 (2011) 2262-2275.

\bibitem{yu2011} W. Yu, G. Chen, M. Cao, Consensus in directed networks of agents with nonlinear dynamics, IEEE Trans. on Automat. Contr. 56 (2011) 1436-1441.

\bibitem{karl2012} M. Andreasson, D. Dimarogonas, K. Johansson, Undamped nonlinear consensus using integral lyapunov function, in Proc. American Control Conference (ACC) (2012) 6644 - 6649.

\bibitem{xie2014} D.M. Xie, Q.L. Liu, L.F. Lv, S.Y. Li, Necessary and sufficient condition for the group consensus of multi-agent systems, Appl. Math. Comput. 243 (2014), 870–878.

\bibitem{WangChenMa2014} J. Wang, K. Chen and Q. Ma, Adaptive Leader-Following Consensus of Multi-Agent Systems
with Unknown Nonlinear Dynamics, Entropy 16 (2014) 5020-5031.


\bibitem{wang2015} Z. Wang, L. Wang, A. Szolnoki, M. Perc, Evolutionary games on multilayer networks: A colloquium, Phys. Condensed Matter 88 (2015) 60270-60277. 

\bibitem{zhaodw2015} D.W. Zhao, L.H. Wang, L.J. Xu, Z. Wang, Finding another yourself in multiplex networks, Applied Math. Comp. 266 (2015) 599-604.

\bibitem{zhao2015} L. Zhao, Y.M. Jia, Finite-time consensus for second-order stochastic multi-agent systems with nonlinear dynamics, Appl. Math. Comput. 270 (2015) 278–290.

\bibitem{LiDuanChenHuang2010} Z. Li, Z. Duan, G. Chen and L. Huang, Consensus of Multiagent Systems and
Synchronization of Complex Networks: A Unified Viewpoint, IEEE Trans. Circuit and Sys. 57 (2010) 213-224.

\bibitem{ZhuYuan2014} J. Zhu and L. Yuan, Consensus of high-order multi-agent systems with switching topologies, Linear Algebra and its Applications 443 (2014) 105-119.


\bibitem{MovricLewis2014} K. H. Movric, F. L. Lewis, F. L., Cooperative Optimal Control for Multi-Agent Systems
on Directed Graph Topologies, IEEE Trans. on Automat. Contr. 59 (2014) 769-774 . 

\bibitem{BaussoGiarrePesenti2006} D. Bauso, L. Giarre, L. R. Pisanti, Non-linear protocols for optimal distributed consensus in networks of dynamic agents, Systems \& Control Letters 55 (2006) 918-928.

\bibitem{QuChunyuWang2007} Z. Qu, C. J, J. Wang, Nonlinear Cooperative Control for Consensus of Nonlinear and
Heterogeneous Systems, in Proc. 46th IEEE Conf. Decision and Control (2007) 2301-2308.

\bibitem{LiuXieRenWang2013} K. Liu, G. Xie, W. Ren and L. Wang, Consensus for multi-agent systems with inherent nonlinear dynamics under
directed topologies, Systems \& Control Letters 62 (2013) 152-162.

\bibitem{XuCaoYuLu2013} W. Xu, J. Cao, W. Yu and J. Lu, Leader-following consensus of non-linear multi-agent
systems with jointly connected topology, IET Control Theory Appl. 8 (2014) 432-440.

\bibitem{dunbar2005} W. Dunbar, A Distributed Receding Horizon Control Algorithm for Dynamically Coupled Nonlinear Systems, in Proc. 44th IEEE Conf. Decision and Control, and the European Control Conference (2005) 6673-6679.

\bibitem{LiShi2013} H. Li, Y. Shi, Distributed model predictive control of constrained nonlinear
systems with communication delays, Systems \& Control Letters 62 (2013) 19-26.

\bibitem{shi2014} H. Li, Y. Shi, Robust distributed model predictive control of constrained continuous-time nonlinear systems: a robustness constraint approach, IEEE Trans. on Automat. Contr. 59 (2014) 1673-1678.

\bibitem{QiuDuan2014} H. Qiu, H. Duan, Receding horizon control for multiple UAV formation flight
based on modified brain storm optimization,  Nonlinear Dyn. 78 (2014) 1973-1988.

\bibitem{bryson1975} A. Bryson, Y. Ho, Applied optimal control, Hemisphere, New York, 1975  Secs. 2.3 and 6.3.

\bibitem{kabamba1987} P. Kabamas, R. Longman, S. Jian-Guo, A homotopy approach to the feedback stablilization of linear systems, J. Guidance, Control, Dyn. 10 (1987) 422-432.

\bibitem{ohtsuka1994} T. Ohtsuka, H. Fujii HA, Stabilized continuation method for solving optimal control problems, J. Guidance, Control, Dyn 17 (1994) 950-957.

\bibitem{ohtsuka1997} T. Ohtsuka, H. Fujii, Real-Time Optimization Algorithm for Nonlinear Receding-Horizon Control,  Automatica 33 (1997) 1147-1154.

\bibitem{ohtsuka1998} T. Ohtsuka, Time-variant receding-horizon control of nonlinear systems, J. Guidance, Control, Dyn. 21 (1998) 174-176.

\bibitem{antoniou2007} A. Antoniou, W. Lu, Practical Optimization: Algorithm and Engineering Applications, Springer, 2007, sec. 10.8.
    
\bibitem{jadbabaie2005} A. Jadbabaie, J. Hauser, On the stability of receding horizon control with a general terminal cost', IEEE Trans. Autom. Contr. 50 (2005) 674-678.

\bibitem{lorenz1963} E. Lorenz, Deterministic Nonperiodic Flow,  J. Atmos. Sci. 20 (1963) 130-141.

\bibitem{lu_chen2002} J. L\"{u} and G. Chen, A New Chaotic Attractor Coined, Int. J. Bifurcation \& Chaos 12 (2002) 659-661.

\bibitem{chen1999} G. Chen, T. Ueta, Yet another chaotic attractor, Int. J. Bifurcation \& Chaos 9 (1999) 1465-1466.


\end{thebibliography}
\end{document}